\documentclass[12pt]{article}
\usepackage{bbm}
\usepackage{amsfonts}
\usepackage{mathrsfs}

\usepackage{bbding}
\usepackage{amssymb}
\usepackage{amsmath}
\usepackage{graphicx}
\usepackage[numbers,sort&compress]{natbib}

\newcommand{\bea}{\begin{eqnarray}}
\newcommand{\eea}{\end{eqnarray}}
\newcommand{\Bea}{\begin{eqnarray*}}
\newcommand{\Eea}{\end{eqnarray*}}
\newcommand{\ba}{\begin{array}}
\newcommand{\ea}{\end{array}}
\newcommand{\bt}{\begin{tabular}}
\newcommand{\et}{\end{tabular}}
\newcommand{\btb}{\begin{table}}
\newcommand{\etb}{\end{table}}
\newcommand{\bc}{\begin{center}}
\newcommand{\ec}{\end{center}}
\newcommand{\beq}{\begin{equation}}
\newcommand{\eeq}{\end{equation}}
\newcommand{\bpr}{\begin{proof}}
\newcommand{\epr}{\end{proof}}

\newtheorem{lem}{Lemma}
\newtheorem{thm}{Theorem}
\newtheorem{assum}{Assumption}
\newtheorem{defn}{Definition}

\begin{document}

\title{Fully Coupled Forward-backward Stochastic Differential Equations on Markov Chains }\footnotetext[1]{Xiao's research is supported by NNSF project(11301309) of China, NSFC(11001155) and Shandong Provincial Scientific Research Foundation for Excellent Young Scientists(BS2013SF003)
.}\footnotetext[2]{Corresponding author: Xinling
Xiao}\footnotetext[3] {E-mail address:jsl\symbol{64}sdu.edu.cn; liuhaodong\symbol{64}mail.sdu.edu.cn;
xinlingxiao\symbol{64}mail.sdu.edu.cn.}

\author{Shaolin Ji$^{1}$ Haodong Liu$^{1}$ Xinling Xiao$^{2}$   \\[3mm]
{\small $^1$ Qilu Institute of Finance, Shandong University,
 Jinan 250100, P. R. China}\\
{\small $^2$School of Mathematical Sciences, Shandong Normal
University, Jinan 250014, P. R. China}}
\date{}
\maketitle

\begin{abstract}
We define fully coupled forward-backward stochastic differential equations on spaces related to continuous time,
finite state Markov Chains. Existence and uniqueness results of the fully coupled forward-backward stochastic differential
equations on Markov Chains are obtained.
\end{abstract}

\newpage
\section{Introduction}

Since the first introduction by Pardoux and Peng \cite{Peng01} in 1990, the theory of nonlinear backward stochastic differential Equations(BSDEs, for short) driven by a Brownian motion has been intensively researched by many researchers and has been achieved abundant theoretical results.  Now the theory of BSDEs has risen to be a powerful part of the stochastic analyst's tool. It also has many important applications, namely in stochastic control, stochastic differential game, finance and the theory of partial differential equations(PDEs, for short).

In the classic BSDEs theory, we consider the Brownian motion as the driver, but Brownian motion is a kind of very idealized stochastic model which restricts greatly the applications  of the classic BSDEs. There are many theories of BSDEs with jump processes recently. Tang and Li \cite{Tang jump} first discussed BSDEs driven by Brownian motion and Poisson Process; Nualart, Schoutens \cite{Levy1} gave BSDEs driven by Brownian motion and Levy Process. Furthermore there are results on BSDEs that other process instead of the Brown motion in diffusion term, such as Ouknine \cite{P-5} researched BSDE driven by Poisson random measure. In 2008, Cohen \cite{Cohen1} studied  BSDE driven by a continuous time, finite state Markov Chain. After then, many results such as comparison theorem about this kind of BSDE,  nonlinear expected results(\cite{Cohen2}, \cite{Cohen3}) and so on were achieved.

Along with the rapid development of the BSDE theory, the theory of fully coupled forward backward stochastic differential equations(FBSDEs, for short) which is closely related to BSDEs has been developing very rapidly. Fully coupled FBSDEs with Brownian motion can be encountered in the optimization problem when applying stochastic maximum principle and mathematical finance considering large investor in security market. As we know now, to get results on the existence and the uniqueness of fully coupled FBSDE's solutions, there are mainly three methods: method of contraction Mapping(\cite{Antonelli}, \cite{Pardoux Tang 1999 Viscosity}), four-step scheme(\cite{4 steps}, \cite{Ma-Yong FBSDE}), method of continuation(\cite{HuPeng}, \cite{Peng Wu 1999}). For more contents on fully coupled FBSDEs, the reader is refereed to Yong \cite{Yong FBSDE2}, Wu \cite{Wu FBSDE 2} or Ma, Wu, Zhang and Zhang \cite{Ma-Zhang FBSDE} and the references therein. Wu (\cite{Wu P-1}, \cite{Wu P-2})also studied FBSDEs driven by Brownian motion and Poisson Process. Li \cite{Li Na} studied FBSDEs driven by Brownian motion and Levy Process.

In this paper, we study fully coupled FBSDEs driven by a Martingale which is generated by a continuous time, finite state Markov Chain  when SDE and BSDE have different dimensions. Inspired by Peng and Wu \cite{Peng Wu 1999}, we also introduced a $m\times n$ full rank matrix $G$ to overcome the problem caused by the different dimensions of SDE and BSDE. Using the method of continuation, the $It\hat{o}$ product rule of Semimartingales and the fixed point principle, with the help of the theory of BSDE driven by a continuous time, finite state Markov Chain, We obtain the existence and uniqueness results of the FBSDEs on Markov Chains. Because of the property of the martingale generated by the finite state Markov Chain is different from the property of the Brownian motion, so the form of the monotone assumptions we got here is different from Peng and Wu \cite{Peng Wu 1999}.

In Section 2, we consider formulation of the problem; in Section 3, we give preliminary concerns; in Section 4 we study the existence and uniqueness results of the  fully coupled FBSDEs on Markov chains under some monotone assumptions; in Section 5 we give proofs of Lemmas in Section 4; in Section 6 we discuss the existence and uniqueness results of the  fully coupled FBSDEs on Markov chains under another set of monotone assumptions; in Section 7 we give proofs of Lemmas in Section 6.

\section{Formulation of the problem}

Assume $m=\{m_{t},t\in[0,T]\}$ is a continuous time, finite state
Markov chain and it takes values in unit vectors $e_{i}$ in $R^{d}$,
where $d$ is the number of states of the chain.

We consider stochastic processes defined on the filtered probability
space $(\Omega,\mathscr{F},\mathscr{F}_{t},\mathbb{P})$, where
$\{\mathscr{F}_{t}\}$ is the completed natural filtration generated
by the $\sigma$-fields $\mathscr{F}_{t}=\sigma(\{m_{s},s\leq t
\},F\in\mathscr{F}_{T}:\mathbb{P}(F)=0)$, and
$\mathscr{F}=\mathscr{F}_{T}$. Observe that $m$ is right-continuous,
this filtration is right-continuous. If $\mathcal {A}_{t}$ denotes
the rate matrix for $m$ at time $t$, then this chain has the
following representation: $$m_{t}=m_{0}+\int_{0}^{t}\mathcal
{A}_{s}m_{s}ds+M_{t}$$

where $M_{t}$ is a martingale (see \cite{Hidden Markov}).

We consider the following forward-backward stochastic differential
equations:
\begin{equation}\label{fa1}
  \left\{
  \begin{aligned}
    X_{t} &= x + \int_{0}^{t}b(s,X_{s},Y_{s},Z_{s})ds+
    \int_{0}^{t}\sigma(s-,X_{s-},Y_{s-},Z_{s-})dM_{s},  \\
    Y_{t} &= \Phi(X_{T}) + \int_{t}^{T}f(s,X_{s},Y_{s},Z_{s})ds-
    \int_{t}^{T}Z_{s-}dM_{s},  \\
\end{aligned}
 \right.
\end{equation}
where $X,Y,Z$ take values in
$\mathbb{R}^{n},\mathbb{R}^{m},\mathbb{R}^{m\times d}$, $T>0$ is an
arbitrarily fixed number and $b,\sigma,f,\Phi$ are functions with
following appropriate dimensions:
$$b: \Omega\times [0,T]\times \mathbb{R}^{n}\times \mathbb{R}^{m}\times \mathbb{R}^{m\times d}\longrightarrow\mathbb{R}^{n}$$
$$\sigma: \Omega\times [0,T]\times \mathbb{R}^{n}\times \mathbb{R}^{m}\times \mathbb{R}^{m\times d}\longrightarrow\mathbb{R}^{n\times d}$$
$$f: \Omega\times [0,T]\times \mathbb{R}^{n}\times \mathbb{R}^{m}\times \mathbb{R}^{m\times d}\longrightarrow\mathbb{R}^{m}$$
$$\Phi: \Omega\times \mathbb{R}^{n}\longrightarrow \mathbb{R}^{m}.$$

We seek a $\mathscr{F}_{t}$-adapted triple $(X,Y,Z)$ such that it
satisfies the above forward-backward stochastic differential
equations, on$[0,T]$ , $\mathbb{P}$-almost surely. That is, our aim
is to find the $\mathscr{F}_{t}$-adapted solution of the the above
forward-backward stochastic differential equations.

Note that, as it is only the $s-left$ limit $\sigma(s-,X_{s-},Y_{s-},Z_{s-})$ which enters into Eqs(\ref{fa1}),
there is no loss of generality to assume that $\sigma(s,X_{s},Y_{s},Z_{s})$ is left continuous in $s$ for each $w$ and $X,Y.$
Note also that as $M$ is a semimartingale, $X,Y$ is $c\grave{a}dl\grave{a}g$ and adapted. We suppose the
existence of the left limits of $Z$. So $Z$ must have at most a countable number of discontinuities, and
then it must be left-continuous at each $t$ except possibly on a $dt$-null set. Hence, if $Z_{s}$ satisfies
Eqs(\ref{fa1}), then so does $Z_{s-}$,

\begin{equation*}
  \left\{
  \begin{aligned}
    X_{t} &= x + \int_{0}^{t}b(s,X_{s},Y_{s},Z_{s-})ds+
    \int_{0}^{t}\sigma(s-,X_{s-},Y_{s-},Z_{s-})dM_{s}, \\
    Y_{t} &=  \Phi(X_{T}) + \int_{t}^{T}f(s,X_{s},Y_{s},Z_{s-})ds-
    \int_{t}^{T}Z_{s-}dM_{s}.  \\
\end{aligned}
 \right.
\end{equation*}

Writing $Z_{t}^{**}:=Z_{t-}$ we have a left-continuous process
$Z^{**}$ which also satisfy the desired equations, and therefore the
writing of the left limits $Z_{t-}$ is unnecessary as we simply
assume our solution is left-continuous.

Given there arguments, we rewrite Eqs(\ref{fa1}) as

\begin{equation*}
  \left\{
  \begin{aligned}
    X_{t} &= x + \int_{0}^{t}b(s,X_{s},Y_{s},Z_{s})ds+
    \int_{0}^{t}\sigma(s,X_{s-},Y_{s-},Z_{s})dM_{s}, \\
    Y_{t} &= \Phi(X_{T}) + \int_{t}^{T}f(s,X_{s},Y_{s},Z_{s})ds-
    \int_{t}^{T}Z_{s}dM_{s}.  \\
\end{aligned}
 \right.
\end{equation*}

As a side note, observe that Eqs(\ref{fa1}) is equivalent to
\begin{equation} \label{fa2}
  \left\{
  \begin{aligned}
    X_{t} &= x + \int_{0}^{t}b^{**}(s,X_{s},Y_{s},Z_{s})ds+
    \int_{0}^{t}\sigma(s,X_{s-},Y_{s-},Z_{s})dm_{s}, \\
    Y_{t} &= \Phi(X_{T}) + \int_{t}^{T}f^{**}(s,X_{s},Y_{s},Z_{s})ds-
    \int_{t}^{T}Z_{s}dm_{s},  \\
\end{aligned}
 \right.
\end{equation}

Where
$b^{**}(s,X_{s},Y_{s},Z_{s})=b(s,X_{s},Y_{s},Z_{s})-\sigma(s,X_{s},Y_{s},Z_{s})\mathcal
{A}_{s}m_{s}(\omega)$ and

$f^{**}(s,X_{s},Y_{s},Z_{s})=f(s,X_{s},Y_{s},Z_{s})+Z_{s}\mathcal
{A}_{s}m_{s}(\omega)$.

Obviously, Eqs(\ref{fa2}) is  driven by Markov chain $m=\{m_{t},t\in[0,T]\}$.

\section{preliminary concerns  }

\bigskip
{\bf (1) preliminary notations }
\bigskip

We know that the optional quadratic variation of $M_{t}$ is given
by the matrix process $$[M,M]_{t}={\sum}_{s=0}^{t} \bigtriangleup M_{s} \triangle M_{s}^{*},$$
where $``^{*}"$ means transpose.

Observing that $\mathcal {A}$ is the rate matrix of the Markov chain
$m,$ the predictable quadratic variation is
$$\langle M,M\rangle _{t}=\int_{0}^{t}[diag(\mathcal {A}_{s}m_{s})-diag(m_{s})\mathcal {A}_{s}^{*}-\mathcal {A}_{s}diag(m_{s})]ds.$$
Then it can be seen through considering
$$diag(m_{t})=diag(m_{0})+\int_{0}^{t}diag(\mathcal {A}_{s}m_{s})ds+\int_{0}^{t}diag(dM_{s})$$

and
$$diag(m_{t})=m_{t}m_{t}^{*}~~~~~~~~~~~~~~~~~~~~~~~~~~~~~~~~~~~~~~~~~~~~~~~~~~~~~~~~~~~$$
    $$~~~~~~~~~~~~~~~~~~~=m_{0}m_{0}^{*}+\int_{0}^{t}m_{s}m_{s}^{*}\mathcal {A}_{s}^{*}ds+\int_{0}^{t}m_{s-}dM_{s}^{*}+\int_{0}^{t}\mathcal {A}_{s}m_{s}m_{s}^{*}ds$$
$$~~~~~~~~~~~~~+\int_{0}^{t}dM_{s}m_{s-}^{*}+{\sum}_{s=0}^{t} \bigtriangleup
M_{s} \triangle M_{s}^{*}.$$

Equating these two we can get
$$[M,M]_{t}=L_{t}+\int_{0}^{t}[diag(\mathcal
{A}_{s}m_{s})-diag(m_{s})\mathcal {A}_{s}^{*}-\mathcal
{A}_{s}diag(m_{s})]ds,$$ where $L$ is some martingale. This in turn
suggests that
$$\langle M,M\rangle _{t}=\int_{0}^{t}[diag(\mathcal
{A}_{s}m_{s})-diag(m_{s})\mathcal {A}_{s}^{*}-\mathcal
{A}_{s}diag(m_{s})]ds.$$

 We will define the following notations:

 $(,)$denotes the usual
 inner product in $\mathbb{R}^{n}$ or $\mathbb{R}^{m};$  we use the usual Euclidean norm
 in $\mathbb{R}^{n}$ and $\mathbb{R}^{m};$ and for $z\in \mathbb{R}^{m\times d}$, we define $|z|=\{tr(zz^{*})\}^{\frac{1}{2}}$.

 For $z^{1}\in\mathbb{R}^{m\times d}, z^{2}\in\mathbb{R}^{m\times d},$
 $$((z^{1}, z^{2}))=tr(z^{1}(z^{2})^{*}),$$
 and for $u^{1}=(x^{1},y^{1},z^{1})\in\mathbb{R}^{n}\times \mathbb{R}^{m}\times \mathbb{R}^{m\times d}
 $, $u^{2}=(x^{2},y^{2},z^{2})\in\mathbb{R}^{n}\times \mathbb{R}^{m}\times \mathbb{R}^{m\times d}
 $, $$[u^{1},u^{2}]=(x^{1}, x^{2})+(y^{1}, y^{2} )+((z^{1}, z^{2}))$$

We are given an $m\times n$ full-rank matrix $G$.
For $u=(x,y,z)\in\mathbb{R}^{n}\times \mathbb{R}^{m}\times
\mathbb{R}^{m\times d},$ let $$F(t,u)=(-G^{*}f(t,u), Gb(t,u), 0),$$
$$H(t,u)=(0, 0, G\sigma(t,u))$$
where $G\sigma=(G\sigma_{1} \cdots G\sigma_{d}).$

Furthermore we define
$$\langle C,D\rangle _{V}=Tr(C[diag(\mathcal
{A}_{s}v)-diag(v)\mathcal {A}_{s}^{*}-\mathcal
{A}_{s}diag(v)]D^{*}),$$
$$\|C\|_{v}^{2}=\langle C,C\rangle _{v},$$
where $v$ is a basis vector in $\mathbb{R}^{d}.$

\bigskip
{\bf (2) BSDEs on Markov chains }
\bigskip

We first discuss the existence and uniqueness of a solution to the
following BSDEs on Markov chains. We consider the equations of the
form
$$Y_{t} = \Phi(X_{T}) - \int_{t}^{T}F_{1}(\omega,s,Y_{s},Z_{s})ds-
    \int_{t}^{T}[F_{2}(s,Y_{s-})+Z_{s}]dM_{s}$$
for functions $F_{1}:\Omega \times [0,T]\times \mathbb{R}^{n}\times
\mathbb{R}^{n \times d}$ and $F_{2}:\Omega \times [0,T]\times
\mathbb{R}^{n}.$ These functions are assumed to be progressively
measurable, $i.e.$ $F_{1}(\cdot,s,Y_{s},Z_{s})$ and
$F_{2}(\cdot,s,Y_{s})$ are $\mathscr{F}_{t}$ measurable for all
$t\in [0,T] $.

From Cohen \cite{Cohen1}, we have the following result about martingale
representation:

\begin{lem}\label{lemma1}
 Any $R^{d}$ valued martingale $L$ defined on
$(\Omega,\mathscr{F}_{t},\mathbb{P})$ can be represented as a
stochastic ( in this case Stieltjes) integral with respect to the
martingale process $M$, up to equality $\mathbb{P}-a.s.$. This
representation is unique up to a $d\langle M,M \rangle_{t}\times
\mathbb{P}-null$ set. That is
$$L_{t}=L_{0}+\int_{0}^{t}Z_{s}dM_{s},$$  where $Z_{s}$ is a
predictable $\mathbb{R}^{n\times d}$ valued matrix process.
\end{lem}

Cohen \cite{Cohen1} gave the existence and uniqueness of a solution to the above
BSDEs.

 \begin{thm}\label{theorem1}
  Assume Lipschitz continuity on the generators
 $F_{2}$ and $F_{2}$. we shall require there to exist $c\in\mathbb{R}$  such that
 for all $s\in [0,T]$
 $$E|F_{1}(s,Y_{s}^{1},Z_{s}^{1})-F_{1}(s,Y_{s}^{2},Z_{s}^{2})|^{2}\leq c^{2}E|Y_{s}^{1}-Y_{s}^{2}|^{2}+c^{2}E\|Z_{s}^{1}-Z_{s}^{2}\|_{M_{s}}^{2}$$
 $$E\|F_{2}(s,Y_{s-}^{1})-F_{2}(s,Y_{s-}^{2})\|_{M_{s}}^{2}\leq c^{2}E|Y_{s}^{1}-Y_{s}^{2}|^{2}.$$
 Under the above Lipschitz condition,  the above  equation has
 at most one solution up to indistinguishability for $Y$ and
 equality $a.s. d\langle M,M\rangle _{t}\times \mathbb{P}$
 for $Y$.
\end{thm}

\section{Existence and uniqueness for above forward-backward equations  }

We give two definitions as following:
\begin{defn}
We denote by $\mathcal
{M}^{2}(0,T;\mathbb{R}^{n})$ the set of all ${R}^{n}$-valued
$\mathscr{F}_{t}$-adapted processes such that
$$E\int_{0}^{T}|v(s)|^{2}ds<+\infty.$$
\end{defn}

\begin{defn}
 A triple of processes $(X,Y,Z): \Omega\times
[0,T] \rightarrow \mathbb{R}^{n}\times \mathbb{R}^{m}\times
\mathbb{R}^{m\times d}$ is called an adapted solution of the Eqs
$(1)$ and $(2)$, if $(X,Y,Z)\in\mathcal
{M}^{2}(0,T;\mathbb{R}^{n}\times \mathbb{R}^{m}\times
\mathbb{R}^{m\times d})$ ,and it satisfies Eqs $(1)$ and $(2)$
$\mathbb{P}$-almost surely.
\end{defn}

The solution's adaptedness  allows us to put Eqs(\ref{fa1}) in
a differential form below:

\begin{equation*}
  \left\{
  \begin{aligned}
    dX_{t} &= b(t,X_{t},Y_{t},Z_{t})dt+\sigma(t,X_{t-},Y_{t-},Z_{t})dM_{t}, \\
    -dY_{t} &= f(t,X_{t},Y_{t},Z_{t})dt-Z_{t}dM_{t}, \\
    X_{0}  &=x, Y_{T}=\Phi(X_{T}).\\
\end{aligned}
 \right.
\end{equation*}

Now we give the critical assumptions of our paper:

\begin{assum}\label{Assumption1}
For each $u=(x,y,z)\in\mathbb{R}^{n}\times
\mathbb{R}^{m}\times \mathbb{R}^{m\times d}, F(\cdot ,u), H(\cdot
,u)\in\mathcal {M}^{2}(0,T;\mathbb{R}^{n}\times \mathbb{R}^{m}\times
\mathbb{R}^{m\times d})$  and for each $x\in\mathbb{R}^{n},
\Phi(x)\in L^{2}(\Omega,\mathscr{F}_{T};\mathbb{R}^{n});$ and
 $$(1) F(t,u)~ is~ uniformly~ lipschitz~ with ~respect~ to ~u;$$
 $$(2) H(t,u)~ is~ uniformly~ lipschitz~ with~ respect~ to~ u;$$
 $$(3) \Phi(x)~ is~ uniformly~ lipschitz~ with~ respect~ to~ x.~~~$$
 \end{assum}

\begin{assum}\label{Assumption2}
There exists a constant $c_{2}>0,c_{2}^{'}>0$
such that
$$[F(t,u^{1})-F(t,u^{2}), u^{1}-u^{2}] \leq
-c_{2}|G(x^{1}-x^{2})|^{2}-c_{2}^{'}(|G^{*}(y^{1}-y^{2})|^{2}+|G^{*}(z^{1}-z^{2})|^{2}),$$
$$[H(t,u^{1})-H(t,u^{2}), u^{1}-u^{2}] \leq
-c_{2}|G(x^{1}-x^{2})|^{2}-c_{2}^{'}(|G^{*}(y^{1}-y^{2})|^{2}+|G^{*}(z^{1}-z^{2})|^{2}),$$
$$\mathbb{P}-a.s.,a.e.t\in\mathbb{R}^{+}, \forall
u^{1}=(x^{1},y^{1},z^{1}), u^{2}=(x^{2},y^{2},z^{2})
\in\mathbb{R}^{n}\times \mathbb{R}^{m}\times \mathbb{R}^{m\times
d},~~~~~~$$  and $$(\Phi(x^{1})-\Phi(x^{2}), G(x^{1}-x^{2}))\geq
c_{3}|G(x^{1}-x^{2})|^{2}, \forall x^{1}\in\mathbb{R}^{n}, \forall
x^{2}\in \mathbb{R}^{n}.$$

Where
 $c_{2}, c_{2}^{'}$ and $c_{3}$  are given positive constants.
 \end{assum}

\bigskip

 In this section, we will give the main result of our paper.

\begin{thm}\label{theorem2}
 Let Assumption \ref{Assumption1} and Assumption \ref{Assumption2} hold, then
there exists a unique adapted solution $(X,Y,Z)$ for Eqs(\ref{fa1}).
\end{thm}

{\bf Proof of Uniqueness:}
 If $U^{1}=(X^{1},Y^{1},Z^{1}), U^{2}=(X^{2},Y^{2},Z^{2})$  are two adapted solutions of Eqs $(1).$  We
 set $$(\hat{X},\hat{Y},\hat{Z})=(X^{1}-X^{2},Y^{1}-Y^{2},Z^{1}-Z^{2}),$$
\begin{equation*}
  \left\{
  \begin{aligned}
    \hat{b}(t) &= b(t,U_{t}^{1})-b(t,U_{t}^{2}), \\
    \hat{\sigma}(t) &= \sigma(t,U_{t}^{1})-\sigma(t,U_{t}^{2}), \\
    \hat{f}(t) &=f(t,U_{t}^{1})-f(t,U_{t}^{2}).\\
\end{aligned}
 \right.
\end{equation*}

Then we have
\begin{equation*}
  \left\{
  \begin{aligned}
    d\hat{X_{t}} &= \hat{b}(t,X_{t},Y_{t},Z_{t})dt+\hat{\sigma}(t,X_{t-},Y_{t-},Z_{t})dM_{t}, \\
    -d\hat{Y_{t}} &= \hat{f}(t,X_{t},Y_{t},Z_{t})dt-\hat{Z_{t}}dM_{t}. \\
\end{aligned}
 \right.
\end{equation*}

 From Assumption \ref{Assumption1}, it follows
 $$E({\sup}_{t\in[0,T]}|\hat{X}_{t}|^{2})+E({\sup}_{t\in[0,T]}|\hat{Y}_{t}|^{2})<+\infty.$$
Using Stieltjes chain rule for products
$$d(G\hat{X_{s}}, \hat{Y_{s}})=\hat{Y_{s-}}d(G\hat{X_{s}})+(G\hat{X_{s-}})d\hat{Y_{s}}+d(G\hat{X_{s}})d\hat{Y_{s}}$$
where $d(G\hat{X_{s}})d\hat{Y_{s}}=d[G\hat{X_{s}},\hat{Y_{s}}]_{t},$
and hence,taking expectation and evaluating at $t=T,$ by Assumption \ref{Assumption2}, we can get
\begin{eqnarray*}
\lefteqn{ \ \  E(\hat{Y_{T}},G\hat{X_{T}})=E(\Phi(X_{T}^{1})-\Phi(X_{T}^{2}),G(X_{T}^{1}-X_{T}^{2}))}
\\ && =E\int_{0}^{T}[(-G^{*}\hat{f},\hat{X_{s}})+(G\hat{b},\hat{Y_{s}})]ds+E\sum_{s=0}^{T}((G\hat{\sigma},\hat{Z_{s}}))\triangle M_{s}\triangle M_{s}^{*}\\
   && =E\int_{0}^{T}[F(s,U^{1})-F(s,U^{2}),U^{1}-U^{2}]ds\\
   && \ \ +E\sum_{s=0}^{T}[H(s,U^{1})-H(s,U^{2}),U^{1}-U^{2}] \triangle M_{s}\triangle M_{s}^{*}\\
   && =E\int_{0}^{T}[F(s,U^{1})-F(s,U^{2}),U^{1}-U^{2}]ds\\
   && \ \ +E\int_{0}^{T}[H(s,U^{1})-H(s,U^{2}),U^{1}-U^{2}]d[M,M]_{s}\\
   && =E\int_{0}^{T}[F(s,U^{1})-F(s,U^{2}),U^{1}-U^{2}]ds\\
   && \ \  +E\int_{0}^{T}[H(s,U^{1})-H(s,U^{2}),U^{1}-U^{2}]d \langle M,M \rangle_{s}\\
   &&  \ \ \leq -c_{2}E\int_{0}^{T}|G\hat{X}|^{2}ds-c_{2}^{'}E\int_{0}^{T}(|G^{*}\hat{Y}|^{2}+|G^{*}\hat{Z}|^{2})ds~\\
   &&  \ \ -c_{2}E\int_{0}^{T}|G\hat{X}|^{2}d \langle M,M \rangle_{s}-c_{2}^{'}E\int_{0}^{T}(|G^{*}\hat{Y}|^{2}+|G^{*}\hat{Z}|^{2})d \langle M,M \rangle_{s}.
\end{eqnarray*}

 we get then
\begin{eqnarray*}
\lefteqn{\ \ \ \  c_{2}(E\int_{0}^{T}|G\hat{X}|^{2}ds+E\int_{0}^{T}|G\hat{X}|^{2}d \langle M,M \rangle_{s})}
\\ && \ \  +c_{2}^{'}[E\int_{0}^{T}(|G^{*}\hat{Y}|^{2}+|G^{*}\hat{Z}|^{2})ds+E\int_{0}^{T}(|G^{*}\hat{Y}|^{2}+|G^{*}\hat{Z}|^{2})d \langle M,M \rangle_{s}]+c_{3}|G\hat{X}|^{2}\\
   &&  \ \ \leq 0.
\end{eqnarray*}

As $c_{2}, c_{2}^{'}$ and $c_{3}$  are given positive constants, then $|G\hat{X}|^{2}\equiv0, |G^{*}\hat{Y}|^{2}\equiv0, |G^{*}\hat{Z}|^{2}\equiv0 $. So $\hat{X}_{s}^{1}\equiv\hat{X}_{s}^{2}, \hat{Y}_{s}^{1}\equiv\hat{Y}_{s}^{2}, \hat{Z}_{s}^{1}\equiv\hat{Z}_{s}^{2}, ~d \langle M,M
 \rangle_{t}\times P-a.s.$




~~~~~~~~~~~~~~~~~~~~~~~~~~~~~~~~~~~~~~~~~~~~~~~~~~~~~~~~~~~~~~~~~~~~~~~~~~~~~~~~~~~~~~~~~~~~~~~~~~~~~~~~~~~~~~~~~$\square$

 We now consider the following family of FBSDEs parametrized by $l\in [0,1]:$
\begin{equation}\label{fa3}
   \left\{
  \begin{aligned}
   dX_{t}^{l} &=[(1-l)c_{2}^{'}(-G^{*}Y_{t}^{l})+lb(t, U_{t}^{l})+\phi_{t}]dt+[(1-l)c_{2}^{'}(-G^{*}Z_{t}^{l})+l\sigma(t, U_{t}^{l})+\psi_{t}]dM_{t}  \\
   -dY_{t}^{l} &=[(1-l)c_{2}GX_{t}^{l}+lf(t, U_{t}^{l})+\gamma_{t}]dt-Z_{t}^{l}dM_{t} \\
  X_{0}^{l} &=x, Y_{T}^{l}=l\Phi(X_{T}^{l})+(1-l)GX_{T}^{l}+\xi, \\
   \end{aligned}
   \right.
  \end{equation}
where $\phi, \psi$ and $\gamma$ are given processes in $M^{2}(0,T)$ with values in $R^{n}, R^{n\times d}$ and $R^{m},$ resp.
Clearly, when $l=1, \xi\equiv 0,$ the existence of Eqs$(3)$ implies that of Eqs(\ref{fa1}).  When $l=0,$ Eqs(\ref{fa3}) becomes
\begin{equation}\label{fa4}
   \left\{
  \begin{aligned}
   dX_{t}^{0} &=[-c_{2}^{'}G^{*}Y_{t}^{0}+\phi_{t}]dt+[-c_{2}^{'}G^{*}Z_{t}^{0}+\psi_{t}]dM_{t}  \\
   -dY_{t}^{0} &=[c_{2}GX_{t}^{0}+\gamma_{t}]dt-Z_{t}^{0}dM_{t} \\
  X_{0}^{l} &=x, Y_{T}^{0}=GX_{T}^{0}+\xi. \\
   \end{aligned}
   \right.
  \end{equation}

For proving  the existence part of the theorem, we first need the
following two lemmas.

\begin{lem}\label{lemma2}

The following equation has a unique solution:
\begin{equation}\label{fa5}
   \left\{
  \begin{aligned}
   dX_{t} &=[-c_{2}^{'}G^{*}Y_{t}+\phi_{t}]dt+[-c_{2}^{'}G^{*}Z_{t}+\psi_{t}]dM_{t}  \\
   -dY_{t} &=[c_{2}GX_{t}+\gamma_{t}]dt-Z_{t}dM_{t} \\
   X_{0} &=x, Y_{T}=\lambda GX_{T}+\xi. \\
   \end{aligned}
   \right.
  \end{equation}
where $\lambda$ is a nonnegative constant.
\end{lem}

\begin{lem}\label{lemma3}

We assume Assumption \ref{Assumption1} and Assumption \ref{Assumption2}. If for an $l_{0}\in [0,1)$ there exists a solution $(X^{l_{0}},Y^{l_{0}},Z^{l_{0}})$ of Eqs$(3)$
then there exists a positive constant $\delta_{0},$ such that for each $\delta \in [0,\delta_{0}],$ there exists a solution $(X^{l_{0}+\delta},Y^{l_{0}+\delta},Z^{l_{0}+\delta})$ of Eqs$(3)$ for $l=l_{0}+\delta.$
\end{lem}

Now we can give

{\bf Proof of Existence}~ From Lemma \ref{lemma2} we see immediately that, when $\lambda=1$ in Eqs(\ref{fa5}), eqs(\ref{fa3}) for $l=0$(that is Eqs(\ref{fa4})) has a unique solution.  It then follows from Lemma \ref{lemma3} that there exists a positive constant $\delta_{0}$ depending on Lipschitz constants, $c_{2}, c_{2}^{'}, c_{3}$ and $T$ such that, for each $\delta \in [0,\delta_{0}],$ Eqs$(3)$ for $l=l_{0}+\delta$ has a unique solution. We can repeat this process for $N-times$ with $1\leq N \delta_{0}<1+\delta_{0}.$ It then follows that, in particular, Eqs(\ref{fa3}) for $l=1$ with $\xi \equiv 0$ has a unique solution. The proof is complete.

 ~~~~~~~~~~~~~~~~~~~~~~~~~~~~~~~~~~~~~~~~~~~~~~~~~~~~~~~~~~~~~~~~~~~~~~~~~~~~~~~~~~~~~~~~~~~~~~~~~~~~~~~~~~~~~~~~~$\square$

\section{Proof of Lemma \ref{lemma2}  and Lemma \ref{lemma3}  }

{\bf Proof of Lemma \ref{lemma2}: }

We observe that the matrix $G$ is of full rank.  The proof of the existence for equation $(5)$ will be divided into two cases: $n\leq m$ and $n>m.$

For the first case, the matrix $G^{*}G$ is of full rank.  We set
  \begin{equation*}
  \left(
   \begin{array}{c}
    X^{'} \\
    Y^{'} \\
   Z^{'}
  \end{array}
  \right)
   = \left (
   \begin{array}{c}
     X \\
    G^{*}Y \\
   G^{*}Z
\end{array}
 \right ),
  \left(\begin{array}{c}
    Y^{''} \\
    Z^{''}
 \end{array}\right)=\left(
   \begin{array}{c}
   (I_{m}-G(G^{*}G)^{-1}G^{*})Y \\
   (I_{m}-G(G^{*}G)^{-1}G^{*})Z
 \end{array}
 \right)
\end{equation*}

Multiplying $G^{*}$ on both sides of the BSDE for $(Y,Z)$ yields
\begin{equation}\label{fa6}
   \left\{
  \begin{aligned}
   dX_{t}^{'} &=[-c_{2}^{'}Y_{t}^{'}+\phi_{t}]dt+[-c_{2}^{'}Z_{t}^{'}+\psi_{t}]dM_{t}  \\
   -dY_{t}^{'} &=[c_{2}G^{*}GX_{t}^{'}+G^{*}\gamma_{t}]dt-Z_{t}^{'}dM_{t} \\
   X_{0}^{'} &=x, Y_{T}^{'}=\lambda G^{*}GX_{T}^{'}+G^{*}\xi. \\
   \end{aligned}
   \right.
  \end{equation}

Similarly, multiplying $(I_{m}-G(G^{*}G)^{-1}G^{*})$ on both sides of the same equation
\begin{equation*}
   \left\{
  \begin{aligned}
   -dY_{t}^{''} &=(I_{m}-G(G^{*}G)^{-1}G^{*})\gamma_{t}dt-Z_{t}^{''}dM_{t}, \\
   Y_{T}^{''} &=(I_{m}-G(G^{*}G)^{-1}G^{*})\xi. \\
   \end{aligned}
   \right.
  \end{equation*}

Obviously the pair $(Y^{''},Z^{''})$ is uniquely determined (by [cohen]). The uniqueness of $(X^{'},Y^{'},Z^{'})$ follows from Theorem \ref{theorem2}.  In order to solve Eqs(\ref{fa6}), we
introduce the following $n\times n$-symmetric matrix-valued ODE, known as the matrix-Riccati equation:
\begin{equation*}
   \left\{
  \begin{aligned}
   -\dot{K}(t) &=-c_{2}^{'}K^{2}+c_{2}G^{*}G, t\in [0,T] \\
      K_{T}    &=\lambda G^{*}G. \\
   \end{aligned}
   \right.
  \end{equation*}

It is well known that this equation has a unique nonnegative solution $K(\cdot) \in C^{1}([0,T);S^{n}).$ Where $S^{n}$ stands for the space of all $n\times n$-symmetric matrices. We then consider the solution $(p,q) \in  M^{2}(0,T;R^{n+n\times d)}$ of the following linear simple BSDE:
\begin{equation*}
   \left\{
  \begin{aligned}
   -dp_{t} &=[-c_{2}^{'}K(t)p_{t}+K(t)\phi_{t}+G^{*}\gamma_{t}]dt+[K(t)\psi_{t}-(I_{n}+c_{2}^{'}K(t))q_{t}]dM_{t}, t\in [0,T]\\
   p_{T} &=G^{*}\xi. \\
   \end{aligned}
   \right.
  \end{equation*}

We now let $X_{t}^{'}$ be the solution of the SDE
\begin{equation*}
   \left\{
  \begin{aligned}
   dX_{t}^{'} &=[-c_{2}^{'}(K(t)X_{t}^{'}+p_{t})+\phi_{t}]dt+[\psi_{t}-c_{2}^{'}q_{t}]dM_{t}, \\
   X_{0}^{'} &=x. \\
   \end{aligned}
   \right.
  \end{equation*}

Then it is easy to check that $(X_{t}^{'},Y_{t}^{'},Z_{t}^{'})=(X_{t}^{'},K(t)X_{t}^{'}+p_{t},q_{t})$ is the solution of equation $(6).$  Once $(X^{'},Y^{'},Z^{'})$ and $(Y^{''},Z^{''})$ are resolved, then the triple $(X,Y,Z)$ is uniquely obtained by
  \begin{equation*}
  \left(
   \begin{array}{c}
    X \\
    Y \\
   Z
  \end{array}
  \right)
   = \left (
   \begin{array}{c}
     X^{'} \\
    G(G^{*}G)^{-1}Y^{'}+Y^{''} \\
   G(G^{*}G)^{-1}Z^{'}+Z^{''}
\end{array}
 \right ).
\end{equation*}

\bigskip

For the second case, the matrix $GG^{*}$ is of full rank.  We set
\begin{equation*}
  \left(
   \begin{array}{c}
    X^{'} \\
    Y^{'} \\
    Z^{'} \\
    X^{''}
  \end{array}
  \right)
   = \left (
   \begin{array}{c}
     GX \\
     Y \\
     Z \\
     (I_{n}-G^{*}(GG^{*})^{-1}G)X
\end{array}
 \right ),
  \end{equation*}

$X^{''}$ is the unique solution of the following linear SDE:
\begin{equation*}
   \left\{
  \begin{aligned}
   dX_{t}^{''} &=(I_{n}-G^{*}(GG^{*})^{-1}G)\phi_{t}dt+(I^{n}-G^{*}(GG^{*})^{-1}G)\varphi_{t}dM_{t}, \\
   X_{0}^{''} &=(I_{n}-G^{*}(GG^{*})^{-1}G)x. \\
   \end{aligned}
   \right.
  \end{equation*}

The triple $(X^{'},Y^{'},Z^{'})$ solves the FBSDE:
\begin{equation}\label{fa7}
   \left\{
  \begin{aligned}
   dX_{t}^{'} &=[-c_{2}^{'}GG^{*}Y_{t}^{'}+G\phi_{t}]dt+[-c_{2}^{'}GG^{*}Z_{t}^{'}+G\psi_{t}]dM_{t},  \\
   -dY_{t}^{'} &=[c_{2}X_{t}^{'}+\gamma_{t}]dt-Z_{t}^{'}dM_{t} \\
   X_{0}^{'} &=Gx, Y_{T}^{'}=\lambda X_{T}^{'}+\xi. \\
   \end{aligned}
   \right.
  \end{equation}

To solve this equation, we introduce  the following $m\times m$-symmetric matrix-valued  matrix-Raccati equation:
\begin{equation*}
   \left\{
  \begin{aligned}
   -\dot{K}(t) &=c_{2}I_{m}-c_{2}^{'}KGG^{*}K, t\in [0,T]\\
      K(T)    &=\lambda I_{m}. \\
   \end{aligned}
   \right.
  \end{equation*}

It is well known that this equation has a unique nonnegative solution $K(\cdot) \in C^{1}([0,T);S^{m}).$ Where $S^{m}$ stands for the space of all $m\times m$-symmetric matrices. We then consider the solution $(p,q) \in  M^{2}(0,T;R^{m+m\times d)}$ of the following linear simple BSDE:
  \begin{equation*}
   \left\{
  \begin{aligned}
   -dp_{t} & =[-c_{2}^{'}K(t)GG^{*}p_{t}+K(t)G\phi_{t}+\gamma_{t}]dt+(K(t)G\psi_{t}-(I_{m}+c_{2}^{'}K(t)GG^{*})q_{t})dM_{t},t\in [0,T], \\
    p_{T}  & =\xi. \\
   \end{aligned}
   \right.
  \end{equation*}

We now let $X_{t}^{'}$ be the solution of the SDE:
\begin{equation*}
   \left\{
  \begin{aligned}
   dX_{t}^{'} &=[-c_{2}^{'}GG^{*}(K(t)X_{t}^{'}+\gamma_{t})+G\phi_{t}]dt+[G\psi_{t}-c_{2}^{'}GG^{*}q_{t}]dM_{t},  \\
    X_{0}^{'} &=Gx. \\
   \end{aligned}
   \right.
  \end{equation*}

Then it is easy to check that$(X_{t}^{'},Y_{t}^{'},Z_{t}^{'})=(X_{t}^{'},K(t)X_{t}^{'}+p_{t},q_{t})$ is the solution of eqs(\ref{fa7}).  Once $(X^{'},X^{''}, Y^{'},Z^{'})$ are resolved, by the definition, the triple $(X,Y,Z)$ is uniquely obtained by
  \begin{equation*}
  \left(
   \begin{array}{c}
    X \\
    Y \\
   Z
  \end{array}
  \right)
   = \left (
   \begin{array}{c}
     G^{*}(GG^{*})^{-1}X^{'}+X^{''}\\
    Y^{'} \\
    Z^{'}
\end{array}
 \right ).
\end{equation*}

The proof is complete.

 ~~~~~~~~~~~~~~~~~~~~~~~~~~~~~~~~~~~~~~~~~~~~~~~~~~~~~~~~~~~~~~~~~~~~~~~~~~~~~~~~~~~~~~~~~~~~~~~~~~~~~~~~~~~~~~~~~$\square$
\bigskip

{\bf Proof of Lemma \ref{lemma3}: }
 Since for each $\phi\in M^{2}(0,T; \mathbb{R}^{n}),
\gamma \in M^{2}(0,T; \mathbb{R}^{m}), \psi\in M^{2}(0,T;
\mathbb{R}^{n\times d}), \xi\in L^{2}(\Omega,\mathscr{F}_{T},
\mathbb{P}),x\in\mathbb{R}^{n}, l_{0}\in [0,1),$ there exists a
unique solution of (3), therefore, for each $U_{s}=(X_{s}, Y_{s},
Z_{s})\in M^{2}(0,T; \mathbb{R}^{n+m+m\times d})$ ~there exists a
unique triple $u_{s}=(x_{s}, y_{s}, z_{s})\in M^{2}(0,T;
\mathbb{R}^{n+m+m\times d})$ satisfying the following FBSDE:
\begin{equation*}
  \left\{
  \begin{aligned}
    dx_{t} &=[(1-l_{0})c_{2}^{'}(-G^{*}y_{t})+l_{0}b(t, u_{t})+\delta b(t, U_{t})+\delta c_{2}^{'}(G^{*}Y_{t})+\phi_{t}]dt \\
           &~~~+[(1-l_{0})c_{2}^{'}(-G^{*}z_{t})+l_{0}\sigma(t, u_{t})+\delta \sigma(t, U_{t})+\delta c_{2}^{'}(G^{*}Z_{t})+\psi_{t}]dM_{t}  \\
   -dy_{t} &=[(1-l_{0})c_{2}Gx_{t}+l_{0}f(t, u_{t})+\delta (-c_{2}GX_{t})+\delta f(t, U_{t})+\gamma_{t}]dt-z_{t}dM_{t} \\
        x_{0} &=x, y_{T}=l_{0}\Phi(x_{T})+(1-l_{0})Gx_{T}+\delta (\Phi(X_{T})-GX_{T})+\xi. \\
\end{aligned}
 \right.
\end{equation*}

We want to prove that the mapping defined by
$$I_{l_{0}+\delta}(U\times X_{T})=u\times x_{T}: M^{2}(0,T; \mathbb{R}^{n+m+m\times d})\times L^{2}(\Omega,\mathscr{F}_{T},
\mathbb{P})~~~~~~~~~~~~~~~~~~~~~~~~$$
$$~~~~~\rightarrow M^{2}(0,T;
\mathbb{R}^{n+m+m\times d})\times L^{2}(\Omega,\mathscr{F}_{T},
\mathbb{P})$$ is contract.

Let $U^{'}\times X_{T}^{'}=(X^{'},Y^{'},Z^{'})\times X_{T}^{'} \in
M^{2}(0,T; \mathbb{R}^{n+m+m\times d})\times
L^{2}(\Omega,\mathscr{F}_{T}, \mathbb{P})$ and $u^{'}\times
x_{T}^{'}=I_{l_{0}+\delta}(U^{'} \times X_{T}^{'}).$
We set
\begin{eqnarray*}
        \hat{U}&=&(\hat{X},\hat{Y},\hat{Z})=(X-X^{'},Y-Y^{'},Z-Z^{'}),\\
        \hat{u}&=&(\hat{x},\hat{y},\hat{z})=(x-x^{'},y-y^{'},z-z^{'}), \\
    \hat{f}_{s}&=&f(s,U_{s})-f(s,U_{s}^{'}), \hat{b}_{s}=b(s,U_{s})-b(s,U_{s}^{'}),\hat{\sigma}_{s}=\sigma(s,U_{s})-\sigma(s,U_{s}^{'}),\\
    \bar{f}_{s}&=&f(s,u_{s})-f(s,u_{s}^{'}), \bar{b}_{s}=b(s,u_{s})-b(s,u_{s}^{'}), \bar{\sigma}_{s}=\sigma(s,u_{s})-\sigma(s,u_{s}^{'}).
\end{eqnarray*}

Then
\begin{equation*}
  \left\{
  \begin{aligned}
    d\hat{x}_{s} &= [(1-l_{0})c_{2}^{'}(-G^{*}\hat{y_{s}})+l_{0}\bar{b}_{s}+\delta\hat{b}_{s}+\delta c_{2}^{'}G^{*}\hat{Y_{s}}]ds \\
                &~~~~~~~~+ [(1-l_{0})c_{2}^{'}(-G^{*}\hat{z_{s}})+l_{0}\bar{\sigma}_{s}+\delta\hat{\sigma}_{s}+\delta c_{2}^{'}G^{*}\hat{Z_{s}}]dM_{s}, \\
    -d\hat{y}_{s} &= [(1-l_{0})c_{2}G\hat{x}_{s}+l_{0}\bar{f}_{s}+\delta(-c_{2}G\hat{X}_{s})+\delta\hat{f}_{s}]ds-\hat{z}_{s}dM_{s}, \\
    \hat{x}_{0}  &=0, \hat{y}_{T}  = l_{0}(\Phi(x_{T})-\Phi(x_{T}^{'}))+(1-l_{0})G\hat{x}_{T}+\delta (\Phi(x_{T})-\Phi(x_{T}^{'})-G\hat{x}_{T}). \\
    \end{aligned}
 \right.
\end{equation*}

Using Stieltjes chain rule for product $(G\hat{x}_{s}, \hat{y}_{s}),$ and hence, taking expectations and evaluating at $t=T$ yields
\begin{eqnarray*}
\lefteqn{l_{0}E(\Phi(x_{T})-\Phi(x_{T}^{'}), G\hat{x}_{T})+(1-l_{0})E(G\hat{x}_{T}, G\hat{x}_{T})+\delta E(\Phi(X_{T})-\Phi(X_{T}^{'})-G\hat{X}_{T}, G\hat{x}_{T})}
\\ &&  =l_{0}E\int_{0}^{T}[F(s,u_{s})-F(s,u_{s}^{'}),\hat{u}_{s}]ds+l_{0}E\int_{0}^{T}[H(s,u_{s})-H(s,u_{s}^{'}),\hat{u}_{s}]d \langle M,M \rangle_{s}\\
   && \ \ \ \    -(1-l_{0})E\int_{0}^{T}[c_{2}|G\hat{x}_{s}|^{2}+c_{2}^{'}|G^{*}\hat{y}_{s}|^{2}]ds-(1-l_{0})E\int_{0}^{T}[c_{2}^{'}|G^{*}\hat{z}_{s}|^{2}]d \langle M,M \rangle_{s}\\
   && \ \ \ \   +\delta E\int_{0}^{T}[c_{2}(G\hat{x}_{s},G\hat{X}_{s})+c_{2}^{'}(G^{*}\hat{y}_{s},G^{*}\hat{Y}_{s})+(\hat{x}_{s},-G^{*}\hat{f}_{s})+(G^{*}\hat{y}_{s}, \hat{b}_{s})]ds\\
   && \ \ \ \   +\delta E\int_{0}^{T}[c_{2}^{'}(G^{*}\hat{Z}_{s},G^{*}\hat{z}_{s})+(\hat{z}_{s},G\hat{\sigma}_{s})]d \langle M,M \rangle_{s}.
\end{eqnarray*}

From Assumption \ref{Assumption1} and Assumption \ref{Assumption2}, $c_{2}>0, c_{3}>0,$
we obtain
\begin{eqnarray*}
\lefteqn{\ \ \ \ (l_{0}c_{3}+(1-l_{0}))E|G\hat{x}_{T}|^{2}+c_{2}E\int_{0}^{T}|G\hat{x}_{s}|^{2}ds+c_{2}^{'}E\int_{0}^{T}|G^{*}\hat{y}_{s}|^{2}ds+l_{0}c_{2}^{'}E\int_{0}^{T}|G^{*}\hat{z}_{s}|^{2}ds} \\ && +l_{0}c_{2}^{'}E\int_{0}^{T}|G^{*}\hat{y}_{s}|^{2}d \langle M,M \rangle_{s}+l_{0}c_{2}E\int_{0}^{T}|G\hat{x}_{s}|^{2} d \langle M, M \rangle_{s}+c_{2}^{'}E\int_{0}^{T}|G^{*}\hat{z}_{s}|^{2} d \langle M, M \rangle_{s} \\
   &&\ \ \ \ \ \leq \delta c_{4}E\int_{0}^{T}(|\hat{U}_{s}|^{2}+|\hat{u}_{s}|^{2})ds+\delta c_{4}E\int_{0}^{T}(|\hat{U}_{s}|^{2}+|\hat{u}_{s}|^{2})d \langle M,M \rangle_{s} \\
   &&\ \ \ \ \ \ \  +\delta c_{4}E|\hat{X}_{T}|^{2}+\delta c_{4}E|\hat{x}_{T}|^{2}.
\end{eqnarray*}

For the difference of the solutions $(\hat{y},\hat{z})=(y-y^{'},z-z^{'}),$ we apply the usual technique to the BSDE, we can get
\begin{eqnarray*}
\lefteqn{\ \ \ \ E\int_{0}^{T}|\hat{y}_{s}|^{2}ds+E\int_{0}^{T}|\hat{z}_{s}|^{2}d \langle M,M \rangle _{s}} \\
 &&\ \  \leq c_{4}\delta E\int_{0}^{T}|\hat{U}_{s}|^{2}ds+c_{4}\delta E|\hat{X}_{T}|^{2}+c_{4}E\int_{0}^{T}|\hat{x}_{s}|^{2}ds+c_{4}E|\hat{x}_{T}|^{2}.
\end{eqnarray*}

$$E\int_{0}^{T}|\hat{y}_{s}|^{2}ds+E\int_{0}^{T}|\hat{z}_{s}|^{2}d \langle M,M \rangle _{s}~~~~~~~~~~~~~~~~$$
$$ \leq c_{4}\delta E\int_{0}^{T}|\hat{U}_{s}|^{2}ds+c_{4}\delta E|\hat{X}_{T}|^{2}+c_{4}E\int_{0}^{T}|\hat{x}_{s}|^{2}ds+c_{4}E|\hat{x}_{T}|^{2}.~~~~~~~~~$$

Similarly, for the difference of the solution $\hat{x}=x-x^{'},$ we apply the usual technique to the forward part and get
\begin{eqnarray*}
\lefteqn{\ \ \ \ \ \sup\limits_{0\leq s \leq T}|\hat{x}_{s}|^{2}\leq c_{4}E\int_{0}^{T}(|\hat{y}_{s}|^{2}+|\hat{z}_{s}|^{2})ds+c_{4} E\int_{0}^{T}|\hat{u}_{s}|^{2}d \langle M,M \rangle _{s}} \\
 &&\ \ \ \ \ \ \ \  +c_{4}\delta E\int_{0}^{T}|\hat{U}_{s}|^{2}ds+c_{4}\delta E\int_{0}^{T}|\hat{U}_{s}|^{2}d \langle M,M \rangle _{s}. \\
 &&  E\int_{0}^{T}|\hat{x}_{s}|^{2}ds\leq c_{4}TE\int_{0}^{T}(|\hat{y}_{s}|^{2}+|\hat{z}_{s}|^{2})ds+c_{4} TE\int_{0}^{T}|\hat{u}_{s}|^{2}d \langle M,M \rangle _{s}\\
 &&\ \ \ \ \ \ \ \  +c_{4}\delta TE\int_{0}^{T}|\hat{U}_{s}|^{2}ds+c_{4}\delta T E\int_{0}^{T}|\hat{U}_{s}|^{2}d \langle M,M \rangle _{s}.
\end{eqnarray*}

Here the constant $c_{4}$ depends on the Lipschitz constants as well as $G, c_{2}, c_{2}^{'}$ and $T.$

Combing the above estimates, it is clear that, we always  have

we can get
\begin{eqnarray*}
\lefteqn{\ \ \ \ E\int_{0}^{T}|\hat{u}_{s}|^{2}ds+E\int_{0}^{T}|\hat{u}_{s}|^{2}d \langle M,M \rangle _{s}+E|\hat{x}_{T}|^{2}} \\
&&\ \   \leq c_{5} \delta \{E\int_{0}^{T}|\hat{U}_{s}|^{2}ds+E\int_{0}^{T}|\hat{U}_{s}|^{2}d \langle M,M \rangle _{s}+E|\hat{X}_{T}|^{2}\},
\end{eqnarray*}

where the constant
$c_{5}$ depending on $c_{4}, G, c_{2}$ and $c_{3}.$ If we choose
$\delta_{0}=\frac{1}{2c_{5}},$ then it is clear that, for each fixed
$\delta \in [0, \delta_{0}],$ the mapping $I_{\alpha+\delta}$ is
contract in the sense that
\begin{eqnarray*}
\lefteqn{\ \ \ \ E\int_{0}^{T}|\hat{u}_{s}|^{2}ds+E\int_{0}^{T}|\hat{u}_{s}|^{2}d \langle M,M \rangle _{s}+E|\hat{x}_{T}|^{2}} \\
&&\ \   \leq \frac{1}{2} \{E\int_{0}^{T}|\hat{U}_{s}|^{2}ds+E\int_{0}^{T}|\hat{U}_{s}|^{2}d \langle M,M \rangle _{s}+E|\hat{X}_{T}|^{2}\},
\end{eqnarray*}

It indicates that this mapping has a unique fixed point
$(u^{l_{0}+\delta})=(x^{l_{0}+\delta}, y^{l_{0}+\delta},
 z^{l_{0}+\delta}),$ which is the solution of equation (3) for
 $l=l_{0}+\delta.$  The proof is complete.

 ~~~~~~~~~~~~~~~~~~~~~~~~~~~~~~~~~~~~~~~~~~~~~~~~~~~~~~~~~~~~~~~~~~~~~~~~~~~~~~~~~~~~~~~~~~~~~~~~~~~~~~~~~~~~~~~~~$\square$

\section{Another existence and uniqueness theorem }

Next, we give the other existence and uniqueness theorem. First, we
give another assumption.

\begin{assum}\label{Assumption3}

There exists a constant $c_{2}>0,c_{2}^{'}>0$
such that
$$[F(t,u^{1})-F(t,u^{2}), u^{1}-u^{2}] \geq
c_{2}|G(x^{1}-x^{2})|^{2}+c_{2}^{'}(|G^{*}(y^{1}-y^{2})|^{2}+|G^{*}(z^{1}-z^{2})|^{2}),$$
$$[H(t,u^{1})-H(t,u^{2}), u^{1}-u^{2}] \geq
c_{2}|G(x^{1}-x^{2})|^{2}+c_{2}^{'}(|G^{*}(y^{1}-y^{2})|^{2}+|G^{*}(z^{1}-z^{2})|^{2}),$$
$$\mathbb{P}-a.s.,a.e.t\in\mathbb{R}^{+}, \forall
u^{1}=(x^{1},y^{1},z^{1}), u^{2}=(x^{2},y^{2},z^{2})
\in\mathbb{R}^{n}\times \mathbb{R}^{m}\times \mathbb{R}^{m\times
d},~~~~~~$$  and $$(\Phi(x^{1})-\Phi(x^{2}), G(x^{1}-x^{2}))\leq
-c_{3}|G(x^{1}-x^{2})|^{2}, \forall x^{1}\in\mathbb{R}^{n}, \forall
x^{2}\in \mathbb{R}^{n}.$$

Where
 $c_{2}, c_{2}^{'}$ and $c_{3}$  are given positive constants.
\end{assum}

\begin{thm}\label{theorem3}
 Let Assumption \ref{Assumption1} and Assumption \ref{Assumption3} hold, then
there exists a unique adapted solution $(X,Y,Z)$ for Eqs(\ref{fa1}).
\end{thm}

 {\bf Proof of Uniqueness:}\
 Using the same procedure as the proof of uniqueness of Theorem \ref{theorem2} and by Assumption \ref{Assumption3} we can get
\begin{eqnarray*}
\lefteqn{ \ \  E(\hat{Y_{T}},G\hat{X_{T}})=E(\Phi(X_{T}^{1})-\Phi(X_{T}^{2}),G(X_{T}^{1}-X_{T}^{2}))}
\\  && =E\int_{0}^{T}[F(s,U^{1})-F(s,U^{2}),U^{1}-U^{2}]ds\\
   && \ \  +E\int_{0}^{T}[H(s,U^{1})-H(s,U^{2}),U^{1}-U^{2}]d \langle M,M \rangle_{s}\\
   &&  \ \ \geq c_{2}E\int_{0}^{T}|G\hat{X}|^{2}ds+c_{2}^{'}E\int_{0}^{T}(|G^{*}\hat{Y}|^{2}+|G^{*}\hat{Z}|^{2})ds\\
   &&  \ \ +c_{2}E\int_{0}^{T}|G\hat{X}|^{2}d \langle M,M \rangle_{s}+c_{2}^{'}E\int_{0}^{T}(|G^{*}\hat{Y}|^{2}+|G^{*}\hat{Z}|^{2})d \langle M,M \rangle_{s}.
\end{eqnarray*}

 We get then
\begin{eqnarray*}
\lefteqn{\ \ \ \  c_{2}(E\int_{0}^{T}|G\hat{X}|^{2}ds+E\int_{0}^{T}|G\hat{X}|^{2}d \langle M,M \rangle_{s})}
\\ && \ \  +c_{2}^{'}[E\int_{0}^{T}(|G^{*}\hat{Y}|^{2}+|G^{*}\hat{Z}|^{2})ds+E\int_{0}^{T}(|G^{*}\hat{Y}|^{2}+|G^{*}\hat{Z}|^{2})d \langle M,M \rangle_{s}]+c_{3}|G\hat{X}|^{2}\\
   &&  \ \ \leq 0.
\end{eqnarray*}

As $c_{2}, c_{2}^{'}$ and $c_{3}$  are given positive constants, then $|G\hat{X}|^{2}\equiv0, |G^{*}\hat{Y}|^{2}\equiv0, |G^{*}\hat{Z}|^{2}\equiv0 $. So $\hat{X}_{s}^{1}\equiv\hat{X}_{s}^{2}, \hat{Y}_{s}^{1}\equiv\hat{Y}_{s}^{2}, \hat{Z}_{s}^{1}\equiv\hat{Z}_{s}^{2}, ~d \langle M,M
 \rangle_{t}\times P-a.s.$



~~~~~~~~~~~~~~~~~~~~~~~~~~~~~~~~~~~~~~~~~~~~~~~~~~~~~~~~~~~~~~~~~~~~~~~~~~~~~~~~~~~~~~~~~~~~~~~~~~~~~~~~~~~~~~~~~$\square$

\bigskip

 We now consider the following family of FBSDEs parametrized by $l\in [0,1]:$
\begin{equation}\label{fa8}
   \left\{
  \begin{aligned}
   dX_{t}^{l} &=[(1-l)c_{2}^{'}(G^{*}Y_{t}^{l})+lb(t, U_{t}^{l})+\phi_{t}]dt+[l\sigma(t, U_{t}^{l})+(1-l)c_{2}^{'}(G^{*}Z_{t}^{l})+\psi_{t}]dM_{t}  \\
   -dY_{t}^{l} &=[-(1-l)c_{2}GX_{t}^{l}+lf(t, U_{t}^{l})+\gamma_{t}]dt-Z_{t}^{l}dM_{t} \\
  X_{0}^{l} &=x, Y_{T}^{l}=l\Phi(X_{T}^{l})+(1-l)GX_{T}^{l}+\xi, \\
   \end{aligned}
   \right.
  \end{equation}

where $\phi, \psi$ and $\gamma$ are given processes in $M^{2}(0,T)$ with values in $R^{n}, R^{n\times d}$ and $R^{m},$ resp.
Clearly, when $l=1, \xi\equiv 0,$ the existence of equation $(8)$ implies that of equation $(1).$ When $l=0,$ equation $(8)$ becomes
\begin{equation}\label{fa9}
   \left\{
  \begin{aligned}
   dX_{t}^{0} &=[c_{2}^{'}G^{*}Y_{t}^{0}+\phi_{t}]dt+[c_{2}^{'}G^{*}Z_{t}^{0}+\psi_{t}]dM_{t}  \\
   -dY_{t}^{0} &=[-c_{2}GX_{t}^{0}+\gamma_{t}]dt-Z_{t}^{0}dM_{t} \\
  X_{0}^{l} &=x, Y_{T}^{0}=GX_{T}^{0}+\xi. \\
   \end{aligned}
   \right.
  \end{equation}


For proving  the existence part of the theorem, we first need the
following two lemmas.

\begin{lem}\label{lemma4}
The following equation has a unique solution:
\begin{equation}\label{fa10}
   \left\{
  \begin{aligned}
   dX_{t} &=[c_{2}^{'}G^{*}Y_{t}+\phi_{t}]dt+[c_{2}^{'}G^{*}Z_{t}+\psi_{t}]dM_{t}  \\
   -dY_{t} &=[-c_{2}GX_{t}+\gamma_{t}]dt-Z_{t}dM_{t} \\
   X_{0} &=x, Y_{T}=\lambda GX_{T}+\xi. \\
   \end{aligned}
   \right.
  \end{equation}
where $\lambda$ is a nonnegative constant.
\end{lem}

\begin{lem}\label{lemma5}
We assume Assumption \ref{Assumption1} and Assumption \ref{Assumption3} hold. If for an $l_{0}\in [0,1)$ there exists a solution $(X^{l_{0}},Y^{l_{0}},Z^{l_{0}})$ of Eqs$(8),$
then there exists a positive constant $\delta_{0},$ such that for each $\delta \in [0,\delta_{0}],$ there exists a solution $(X^{l_{0}+\delta},Y^{l_{0}+\delta},Z^{l_{0}+\delta})$ of Eqs$(8)$ for $l=l_{0}+\delta.$
\end{lem}

Now we can give

{\bf Proof of Existence}~ From Lemma \ref{lemma4}, we see immediately that, when $\lambda=1$ in Eqs(\ref{fa10}),  eqs(\ref{fa8}) for $l=0$(that is Eqs(\ref{fa9})) has a unique solution.  It then follows from Lemma \ref{lemma5} that there exists a positive constant $\delta_{0}$ depending on Lipschitz constants, $c_{2}, c_{2}^{'}, c_{3}$ and $T$ such that, for each $\delta \in [0,\delta_{0}],$ Eqs$(8)$ for $l=l_{0}+\delta$ has a unique solution. We can repeat this process for $N-times$ with $1\leq N \delta_{0}<1+\delta_{0}.$ It then follows that, in particular, FBSDE Eqs(\ref{fa8}) for $l=1$ with $\xi \equiv 0$ has a unique solution. The proof is complete.

 ~~~~~~~~~~~~~~~~~~~~~~~~~~~~~~~~~~~~~~~~~~~~~~~~~~~~~~~~~~~~~~~~~~~~~~~~~~~~~~~~~~~~~~~~~~~~~~~~~~~~~~~~~~~~~~~~~$\square$

\bigskip

\section{Proof of Lemma \ref{lemma4} and Lemma \ref{lemma5} }

{\bf Proof of Lemma \ref{lemma4}:}

We observe that the matrix $G$ is of full rank.  The proof of the existence for equation $(10)$ will be divided into two cases: $n\leq m$ and $n>m.$

For the first case, the matrix $G^{*}G$ is strictly positive.  We set
  \begin{equation*}
  \left(
   \begin{array}{c}
    X^{'} \\
    Y^{'} \\
   Z^{'}
  \end{array}
  \right)
   = \left (
   \begin{array}{c}
     X \\
    G^{*}Y \\
   G^{*}Z
\end{array}
 \right ),
  \left(\begin{array}{c}
    Y^{''} \\
    Z^{''}
 \end{array}\right)=\left(
   \begin{array}{c}
   (I_{m}-G(G^{*}G)^{-1}G^{*})Y \\
   (I_{m}-G(G^{*}G)^{-1}G^{*})Z
 \end{array}
 \right)
\end{equation*}

Multiplying $G^{*}$ on both sides of the BSDE for $(Y,Z)$ yields
\begin{equation}
   \left\{
  \begin{aligned}
   dX_{t}^{'} &=[c_{2}^{'}Y_{t}^{'}+\phi_{t}]dt+[c_{2}^{'}Z_{t}^{'}+\psi_{t}]dM_{t}  \\
   -dY_{t}^{'} &=[-c_{2}G^{*}GX_{t}^{'}+G^{*}\gamma_{t}]dt-Z_{t}^{'}dM_{t} \\
   X_{0}^{'} &=x, Y_{T}^{'}=\lambda G^{*}GX_{T}^{'}+G^{*}\xi. \\
   \end{aligned}
   \right.
  \end{equation}

Similarly, multiplying $(I_{m}-G(G^{*}G)^{-1}G^{*})$ on both sides of the same equation
\begin{equation*}
   \left\{
  \begin{aligned}
   -dY_{t}^{''} &=(I_{m}-G(G^{*}G)^{-1}G^{*})\gamma_{t}dt-Z_{t}^{''}dM_{t}, \\
   Y_{T}^{''} &=(I_{m}-G(G^{*}G)^{-1}G^{*})\xi. \\
   \end{aligned}
   \right.
  \end{equation*}

Obviously the pair $(Y^{''},Z^{''})$ is uniquely determined. The uniqueness of $(X^{'},Y^{'},Z^{'})$ follows from Theorem $7.1.$  In order to solve Eqs$(11),$ we
introduce the following $n\times n$-symmetric matrix-valued ODE, known as the matrix-Raccati equation:
\begin{equation*}
   \left\{
  \begin{aligned}
   -\dot{K}(t) &=-c_{2}^{'}K^{2}+c_{2}G^{*}G, t\in [0,T] \\
      K_{T}    &=\lambda G^{*}G. \\
   \end{aligned}
   \right.
  \end{equation*}

It is well known that this equation has a unique nonnegative solution $K(\cdot) \in C^{1}([0,T);S^{n}).$ Where $S^{n}$ stands for the space of all $n\times n$-symmetric matrices. We then consider the solution $(p,q) \in  M^{2}(0,T;R^{n+n\times d)}$ of the following linear simple BSDE:
\begin{equation*}
   \left\{
  \begin{aligned}
   -dp_{t} &=[-c_{2}^{'}K(t)p_{t}-K(t)\phi_{t}+G^{*}\gamma_{t}]dt+[-K(t)\psi_{t}-(I_{n}+c_{2}^{'}K(t))q_{t}]dM_{t}, t\in [0,T]\\
   p_{T} &=G^{*}\xi. \\
   \end{aligned}
   \right.
  \end{equation*}

We now let $X_{t}^{'}$ be the solution of the SDE
\begin{equation*}
   \left\{
  \begin{aligned}
   dX_{t}^{'} &=[c_{2}^{'}(-K(t)X_{t}^{'}+p_{t})+\phi_{t}]dt+[\psi_{t}+c_{2}^{'}q_{t}]dM_{t}, \\
   X_{0}^{'} &=x. \\
   \end{aligned}
   \right.
  \end{equation*}

Then it is easy to check that $(X_{t}^{'},Y_{t}^{'},Z_{t}^{'})=(X_{t}^{'},-K(t)X_{t}^{'}+p_{t},q_{t})$ is the solution of equation $(11).$  Once $(X^{'},Y^{'},Z^{'})$ and $(Y^{''},Z^{''})$ are resolved, then the triple $(X,Y,Z)$ is uniquely obtained by
  \begin{equation*}
  \left(
   \begin{array}{c}
    X \\
    Y \\
   Z
  \end{array}
  \right)
   = \left (
   \begin{array}{c}
     X^{'} \\
    G(G^{*}G)^{-1}Y^{'}+Y^{''} \\
   G(G^{*}G)^{-1}Z^{'}+Z^{''}
\end{array}
 \right ).
\end{equation*}

\bigskip

For the second case, the matrix $GG^{*}$ is of full rank.  We set
\begin{equation*}
  \left(
   \begin{array}{c}
    X^{'} \\
    Y^{'} \\
    Z^{'} \\
    X^{''}
  \end{array}
  \right)
   = \left (
   \begin{array}{c}
     GX \\
     Y \\
     Z \\
     (I_{n}-G^{*}(GG^{*})^{-1}G)X
\end{array}
 \right ),
  \end{equation*}

$X^{''}$ is the unique solution of the following linear SDE:
\begin{equation*}
   \left\{
  \begin{aligned}
   dX_{t}^{''} &=(I_{n}-G^{*}(GG^{*})^{-1}G)\phi_{t}dt+(I_{n}-G^{*}(GG^{*})^{-1}G)\varphi_{t}dM_{t}, \\
   X_{0}^{''} &=(I_{n}-G^{*}(GG^{*})^{-1}G)x. \\
   \end{aligned}
   \right.
  \end{equation*}

The triple $(X^{'},Y^{'},Z^{'})$ solves the FBSDE:
\begin{equation}
   \left\{
  \begin{aligned}
   dX_{t}^{'} &=[c_{2}^{'}GG^{*}Y_{t}^{'}+G\phi_{t}]dt+[c_{2}^{'}GG^{*}Z_{t}^{'}+G\psi_{t}]dM_{t},  \\
   -dY_{t}^{'} &=[-c_{2}X_{t}^{'}+\gamma_{t}]dt-Z_{t}^{'}dM_{t} \\
   X_{0}^{'} &=Gx, Y_{T}^{'}=\lambda X_{T}^{'}+\xi. \\
   \end{aligned}
   \right.
  \end{equation}

To solve this equation, we introduce  the following $m\times m$-symmetric matrix-valued  matrix-Raccati equation:
\begin{equation*}
   \left\{
  \begin{aligned}
   -\dot{K}(t) &=c_{2}I_{m}-c_{2}^{'}KGG^{*}K, t\in [0,T]\\
      K(T)    &=\lambda I_{m}. \\
   \end{aligned}
   \right.
  \end{equation*}

It is well known that this equation has a unique nonnegative solution $K(\cdot) \in C^{1}([0,T);S^{m}).$ Where $S^{m}$ stands for the space of all $m\times m$-symmetric matrices. We then consider the solution $(p,q) \in  M^{2}(0,T;R^{m+m\times d)}$ of the following linear simple BSDE:
  \begin{equation*}
   \left\{
  \begin{aligned}
   -dp_{t} & =[c_{2}^{'}K(t)GG^{*}p_{t}-K(t)G\phi_{t}+\gamma_{t}]dt+(-K(t)G\psi_{t}-(I_{m}-c_{2}^{'}K(t)GG^{*})q_{t})dM_{t},t\in [0,T], \\
    Y_{T}  & =\xi. \\
   \end{aligned}
   \right.
  \end{equation*}

We now let $X_{t}^{'}$ be the solution of the SDE:
\begin{equation*}
   \left\{
  \begin{aligned}
   dX_{t}^{'} &=[-c_{2}^{'}GG^{*}(-K(t)X_{t}^{'}+\gamma_{t})+G\phi_{t}]dt+[G\psi_{t}-c_{2}^{'}GG^{*}q_{t}]dM_{t},  \\
    X_{0}^{'} &=Gx. \\
   \end{aligned}
   \right.
  \end{equation*}

Then it is easy to check that$(X_{t}^{'},Y_{t}^{'},Z_{t}^{'})=(X_{t}^{'},K(t)X_{t}^{'}+p_{t},q_{t})$ is the solution of equation $(12).$  Once $(X^{'},X^{''}, Y^{'},Z^{'})$ are resolved, by the definition, the triple $(X,Y,Z)$ is uniquely obtained by
  \begin{equation*}
  \left(
   \begin{array}{c}
    X \\
    Y \\
   Z
  \end{array}
  \right)
   = \left (
   \begin{array}{c}
     G^{*}(GG^{*})^{-1}X^{'}+X^{''}\\
    Y^{'} \\
    Z^{'}
\end{array}
 \right ).
\end{equation*}

The proof is complete.

 ~~~~~~~~~~~~~~~~~~~~~~~~~~~~~~~~~~~~~~~~~~~~~~~~~~~~~~~~~~~~~~~~~~~~~~~~~~~~~~~~~~~~~~~~~~~~~~~~~~~~~~~~~~~~~~~~~$\square$

{\bf Proof of Lemma \ref{lemma5}:}

 Since for each $\phi\in M^{2}(0,T; \mathbb{R}^{n}),
\gamma \in M^{2}(0,T; \mathbb{R}^{m}), \psi\in M^{2}(0,T;
\mathbb{R}^{n\times d}), \xi\in L^{2}(\Omega,\mathscr{F}_{T},
\mathbb{P}),x\in\mathbb{R}^{n}, l_{0}\in [0,1),$ there exists a
unique solution of (8), therefore, for each $U_{s}=(X_{s}, Y_{s},
Z_{s})\in M^{2}(0,T; \mathbb{R}^{n+m+m\times d})$ ~there exists a
unique triple $u_{s}=(x_{s}, y_{s}, z_{s})\in M^{2}(0,T;
\mathbb{R}^{n+m+m\times d})$~ satisfying the following FBSDE:
\begin{equation*}
  \left\{
  \begin{aligned}
    dx_{t} &=[(1-l_{0})c_{2}^{'}(G^{*}y_{t})+l_{0}b(t, u_{t})+\delta b(t, U_{t})+\delta (-c_{2}^{'}G^{*}Y_{t})+\phi_{t}]dt \\
           &~~~+[(1-l_{0})c_{2}^{'}(G^{*}z_{t})+l_{0}\sigma(t, u_{t})+\delta \sigma(t, U_{t})-\delta c_{2}^{'}(G^{*}Z_{t})+\psi_{t}]dM_{t}  \\
   -dy_{t} &=[-(1-l_{0})c_{2}Gx_{t}+l_{0}f(t, u_{t})+\delta c_{2}GX_{t}+\delta f(t, U_{t})+\gamma_{t}]dt-z_{t}dM_{t} \\
        x_{0} &=x, y_{T}=l_{0}\Phi(x_{T})+(1-l_{0})Gx_{T}+\delta (\Phi(X_{T})-GX_{T})+\xi. \\
\end{aligned}
 \right.
\end{equation*}

We want to prove that the mapping defined by
$$I_{l_{0}+\delta}(U\times X_{T})=u\times x_{T}: M^{2}(0,T; \mathbb{R}^{n+m+m\times d})\times L^{2}(\Omega,\mathscr{F}_{T},
\mathbb{P})~~~~~~~~~~~~~~~~~~~~~~~~$$
$$~~~~~\rightarrow M^{2}(0,T;
\mathbb{R}^{n+m+m\times d})\times L^{2}(\Omega,\mathscr{F}_{T},
\mathbb{P})$$ is contract.

Let $U^{'}\times X_{T}^{'}=(X^{'},Y^{'},Z^{'})\times X_{T}^{'} \in
M^{2}(0,T; \mathbb{R}^{n+m+m\times d})\times
L^{2}(\Omega,\mathscr{F}_{T}, \mathbb{P})$ and $u^{'}\times
x_{T}^{'}=I_{l_{0}+\delta}(U^{'} \times X_{T}^{'}).$ We set
\begin{eqnarray*}
        \hat{U}&=&(\hat{X},\hat{Y},\hat{Z})=(X-X^{'},Y-Y^{'},Z-Z^{'}),\\
        \hat{u}&=&(\hat{x},\hat{y},\hat{z})=(x-x^{'},y-y^{'},z-z^{'}), \\
    \hat{f}_{s}&=&f(s,U_{s})-f(s,U_{s}^{'}), \hat{b}_{s}=b(s,U_{s})-b(s,U_{s}^{'}),\hat{\sigma}_{s}=\sigma(s,U_{s})-\sigma(s,U_{s}^{'}),\\
    \bar{f}_{s}&=&f(s,u_{s})-f(s,u_{s}^{'}), \bar{b}_{s}=b(s,u_{s})-b(s,u_{s}^{'}), \bar{\sigma}_{s}=\sigma(s,u_{s})-\sigma(s,u_{s}^{'}).
\end{eqnarray*}

Then
\begin{equation*}
  \left\{
  \begin{aligned}
     d\hat{x}_{s} &= [(1-l_{0})c_{2}^{'}(G^{*}\hat{y_{s}})+l_{0}\bar{b}_{s}+\delta\hat{b}_{s}+\delta(- c_{2}^{'}G^{*}\hat{Y_{s}})]ds \\
                &~~~~~~~~+ [(1-l_{0})c_{2}^{'}(G^{*}\hat{z_{s}})+l_{0}\bar{\sigma}_{s}+\delta\hat{\sigma}_{s}-\delta c_{2}^{'}(G^{*}\hat{Z_{s}})]dM_{s}, \\
    -d\hat{y}_{s} &= [-(1-l_{0})c_{2}G\hat{x}_{s}+l_{0}\bar{f}_{s}+\delta c_{2}G\hat{X}_{s}+\delta\hat{f}_{s}]ds-\hat{z}_{s}dM_{s}. \\
    \hat{x}_{0}  &=0, \hat{y}_{T}  = l_{0}(\Phi(x_{T})-\Phi(x_{T}^{'}))+(1-l_{0})G\hat{x}_{T}+\delta (\Phi(x_{T})-\Phi(x_{T}^{'})-G\hat{x}_{T}). \\
    \end{aligned}
 \right.
\end{equation*}

Using Stieltjes chain rule for product $(G\hat{x}_{s}, \hat{y}_{s}),$ and hence, taking expectations and evaluating at $t=T$ yields
\begin{eqnarray*}
\lefteqn{l_{0}E(\Phi(x_{T})-\Phi(x_{T}^{'}), G\hat{x}_{T})+(1-l_{0})E(G\hat{x}_{T}, G\hat{x}_{T})+\delta E(\Phi(X_{T})-\Phi(X_{T}^{'})-G\hat{X}_{T}, G\hat{x}_{T})}
\\ &&  =l_{0}E\int_{0}^{T}[F(s,u_{s})-F(s,u_{s}^{'}),\hat{u}_{s}]ds+l_{0}E\int_{0}^{T}[H(s,u_{s})-H(s,u_{s}^{'}),\hat{u}_{s}]d \langle M,M \rangle_{s}\\
   && \ \ \ \   +(1-l_{0})E\int_{0}^{T}[c_{2}|G\hat{x}_{s}|^{2}+c_{2}^{'}|G^{*}\hat{y}_{s}|^{2}]ds+(1-l_{0})E\int_{0}^{T}[c_{2}^{'}|G^{*}\hat{z}_{s}|^{2}]d \langle M,M \rangle_{s} \\
   && \ \ \ \   +\delta E\int_{0}^{T}[-c_{2}(G\hat{x}_{s},G\hat{X}_{s})-c_{2}^{'}(G^{*}\hat{y}_{s},G^{*}\hat{Y}_{s})+(\hat{x}_{s},-G^{*}\hat{f}_{s})+(G^{*}\hat{y}_{s}, \hat{b}_{s})]ds\\
   && \ \ \ \  +\delta E\int_{0}^{T}[-c_{2}^{'}(G^{*}\hat{Z}_{s},G^{*}\hat{z}_{s})+(\hat{z}_{s},G\hat{\sigma}_{s})]d \langle M,M \rangle_{s}.
\end{eqnarray*}

From Assumption \ref{Assumption1} and Assumption \ref{Assumption3}, $c_{2}>0, c_{3}>0,$
we obtain
\begin{eqnarray*}
\lefteqn{ \ \ (l_{0}c_{3}-(1-l_{0}))E|G\hat{x}_{T}|^{2}+c_{2}E\int_{0}^{T}|G\hat{x}_{s}|^{2}ds+c_{2}^{'}E\int_{0}^{T}|G^{*}\hat{y}_{s}|^{2}ds+l_{0}c_{2}^{'}E\int_{0}^{T}|G^{*}\hat{z}_{s}|^{2}ds} \\
 &&\ +l_{0}c_{2}E\int_{0}^{T}|G\hat{x}_{s}|^{2} d \langle M, M \rangle_{s}+l_{0}c_{2}^{'}E\int_{0}^{T}|G^{*}\hat{y}_{s}|^{2}d \langle M,M \rangle_{s}+c_{2}^{'}E\int_{0}^{T}|G^{*}\hat{z}_{s}|^{2} d \langle M, M \rangle_{s} \\
 &&\ \leq \delta c_{4}E\int_{0}^{T}(|\hat{U}_{s}|^{2}+|\hat{u}_{s}|^{2})ds+\delta c_{4}E\int_{0}^{T}(|\hat{U}_{s}|^{2}+|\hat{u}_{s}|^{2})d \langle M,M \rangle_{s}\\
 && \ \ \ +\delta c_{4}E|\hat{X}_{T}|^{2}+\delta c_{4}E|\hat{x}_{T}|^{2}.
\end{eqnarray*}

For the difference of the solutions $(\hat{y},\hat{z})=(y-y^{'},z-z^{'}),$ we apply the usual technique to the BSDE, we can get
\begin{eqnarray*}
\lefteqn{\ \ \ \ E\int_{0}^{T}|\hat{y}_{s}|^{2}ds+E\int_{0}^{T}|\hat{z}_{s}|^{2}d \langle M,M \rangle _{s}} \\
 &&\ \  \leq c_{4}\delta E\int_{0}^{T}|\hat{U}_{s}|^{2}ds+c_{4}\delta E|\hat{X}_{T}|^{2}+c_{4}E\int_{0}^{T}|\hat{x}_{s}|^{2}ds+c_{4}E|\hat{x}_{T}|^{2}.
\end{eqnarray*}

Similarly, for the difference of the solution $\hat{x}=x-x^{'},$ we apply the usual technique to the forward part and get
\begin{eqnarray*}
\lefteqn{\ \ \ \ \ \sup\limits_{0\leq s \leq T}|\hat{x}_{s}|^{2}\leq c_{4}E\int_{0}^{T}(|\hat{y}_{s}|^{2}+|\hat{z}_{s}|^{2})ds+c_{4} E\int_{0}^{T}|\hat{u}_{s}|^{2}d \langle M,M \rangle _{s}} \\
 &&\ \ \ \ \ \ \ \  +c_{4}\delta E\int_{0}^{T}|\hat{U}_{s}|^{2}ds+c_{4}\delta E\int_{0}^{T}|\hat{U}_{s}|^{2}d \langle M,M \rangle _{s}.
 \end{eqnarray*}
\begin{eqnarray*}
  && \lefteqn{\ \ E\int_{0}^{T}|\hat{x}_{s}|^{2}ds\leq c_{4}TE\int_{0}^{T}(|\hat{y}_{s}|^{2}+|\hat{z}_{s}|^{2})ds+c_{4} TE\int_{0}^{T}|\hat{u}_{s}|^{2}d \langle M,M \rangle _{s}} \\
 &&\ \ \ \ \ \ \ \  +c_{4}\delta TE\int_{0}^{T}|\hat{U}_{s}|^{2}ds+c_{4}\delta T E\int_{0}^{T}|\hat{U}_{s}|^{2}d \langle M,M \rangle _{s}.
\end{eqnarray*}

Here the constant $c_{4}$ depends on the Lipschitz constants as well as $G, c_{2}, c_{2}^{'}$ and $T.$

Combing the above estimates, we can get
\begin{eqnarray*}
\lefteqn{\ \ \ \ E\int_{0}^{T}|\hat{u}_{s}|^{2}ds+E\int_{0}^{T}|\hat{u}_{s}|^{2}d \langle M,M \rangle _{s}+E|\hat{x}_{T}|^{2}} \\
&&\ \   \leq c_{5} \delta \{E\int_{0}^{T}|\hat{U}_{s}|^{2}ds+E\int_{0}^{T}|\hat{U}_{s}|^{2}d \langle M,M \rangle _{s}+E|\hat{X}_{T}|^{2}\},
\end{eqnarray*}

 where the constant
$c_{5}$ depending on $c_{4}, G, c_{2}$ and $c_{3}.$ If we choose
$\delta_{0}=\frac{1}{2c_{5}},$ then it is clear that, for each fixed
$\delta \in [0, \delta_{0}],$ the mapping $I_{\alpha+\delta}$ is
contract in the sense that
\begin{eqnarray*}
\lefteqn{\ \ \ \ E\int_{0}^{T}|\hat{u}_{s}|^{2}ds+E\int_{0}^{T}|\hat{u}_{s}|^{2}d \langle M,M \rangle _{s}+E|\hat{x}_{T}|^{2}} \\
&&\ \   \leq \frac{1}{2} \{E\int_{0}^{T}|\hat{U}_{s}|^{2}ds+E\int_{0}^{T}|\hat{U}_{s}|^{2}d \langle M,M \rangle _{s}+E|\hat{X}_{T}|^{2}\},
\end{eqnarray*}

It indicates that this mapping has a unique fixed point
$(u^{l_{0}+\delta})=(x^{l_{0}+\delta}, y^{l_{0}+\delta},
 z^{l_{0}+\delta}),$ which is the solution of equation (8) for
 $l=l_{0}+\delta.$

 The proof is complete.

 ~~~~~~~~~~~~~~~~~~~~~~~~~~~~~~~~~~~~~~~~~~~~~~~~~~~~~~~~~~~~~~~~~~~~~~~~~~~~~~~~~~~~~~~~~~~~~~~~~~~~~~~~~~~~~~~~~$\square$


\newpage
\begin{thebibliography}{99}
\bibitem{Antonelli}
F. Antonelli, \textit{Backward-forward stochastic differential equations}, The Annals of Applied Probability, 3(1993), 777-793£®

\bibitem{Levy2}
 K. Bahlali, M. Eddahbi, E. Essaky, \textit{BSDE Associated with Levy proceses and application to PDIE}, International Journal of Stochastic Analysis, 16(1)(2003), 1-17.
\bibitem{Cohen1}
S. Cohen, R. Elliott, \textit{Solutions of  Backward Stochastic Differential Equations on Markov Chains}, Communications on Stochastic Analysis, 2(2)(2008), 251-262.
\bibitem{Cohen2}
S. Cohen, R. Elliott, \textit{Comparisons for Backward Stochastic Differential Equations on Markov Chains and related No-Arbitrage Conditions}, The Annals of Applied Probability, 20(1)(2010), 267-311.
\bibitem{Cohen3}
 S. Cohen, Y. Hu, \textit{Ergodic BSDEs driven by Markov Chains}, SIAM Journal on Control and Optimization, 51(5)(2013), 4138-4168.
\bibitem{Hidden Markov}
R. Elliott, L. Aggoun, and J. Moore, \textit{Hidden markov models: Estimation and control}, Springer-Verlag, Berlin-Heidelberg-New York, 1994.
\bibitem{HuPeng}
Y. Hu, S. Peng, \textit{Solutions of forward-backward stochastic differential equations}, Probability Theory and Related Fields, 103(2)(1995), 273-283£®
\bibitem{Li Na}
 N. Li, Z. Wu, \textit{Optimal control problem with delay and Levy processes}, Applied Mathematics-A Journal of Chinese Universities, 29(1)(2014), 67-85.
 \bibitem{Ma-Zhang FBSDE}
 J. Ma, Z. Wu, D. Zhang, J. Zhang, \textit{On wellposedness of forward-backward SDEs-A unified approach}, arxiv.org/abs/math/1110.4658vl, (2011).
 \bibitem{Ma-Yong FBSDE}
 J. Ma, J. Yong, \textit{Forward-backward stochastic diferential equations and their applications}, Lecture Notes in Mathematics, Springer, 1999.
\bibitem{4 steps}
 J. Ma, P. Protter, J. Yong, \textit{Solving forward-backward stochastic differential equations explicitly-a four step scheme}, Probability Theory and Related Fields, 98(3)(1994), 339-359.
 \bibitem{Levy1}
D. Nualart, W. Schoutens, \textit{Backward stochastic differential equations and Feynman-Kac formula for Levy processes, with applications in finance}, Bernoulli, 7(5)(2001):761-776.
 \bibitem{P-5}
 Y. Ouknine, \textit{Reflected Backward Stochastic Differential Equations with Jumps}, Stochastics: An International Journal of Probability and Stochastic Processes, 65(1-2)(1998), 111-125.
 \bibitem{Peng01}
E. Pardoux£¬S. Peng, \textit{Adapted solution of a backward stochastic diferential equations}, System and Control Letters, 14(1990), 55¡ª61£®
\bibitem{Pardoux Tang 1999 Viscosity}
E. Pardoux, S. Tang, \textit{Forward-backward stochastic differential equations and quasilinear parabolic PDEs}, Probability Theory and Related Fields, 114(2)(1999), 123-150.
\bibitem{Peng Wu 1999}
S. Peng, Z. Wu, \textit{Fully coupled forward-backward stochastic differential equations}, SIAM Journal on Control and Optimization, 37(3)(1999), 825-843.
\bibitem{Tang jump}
S. Tang, X. Li, \textit{Necessary conditions for optimal control of stochastic systems with random jumps}, SIAM Journal on Control and Optimization, 32(5)(1994), 1447-1475.
\bibitem{Wu P-1}
Z. Wu, \textit{Forward-Backward stochastic differential equations with Brownian motion and Poisson process}, Acta Mathematicae Applicatae Sinica, 15(4)(1999), 433-443.
\bibitem{Wu P-2}
Z. Wu, \textit{Fully coupled FBSDE with Brownian motion and Poisson process in stopping time duration}, Journal of the Australian Mathematical Society, 74(02)(2003), 249-266.
\bibitem{Wu FBSDE 2}
 Z. Wu, \textit{The comparison theorem of FBSDE}, Statistics and Probability Letters, 44(1)(1999), 1-6.
\bibitem{Yong FBSDE2}
J. Yong, \textit{Forward-backward stochastic differential equations with mixed initial-terminal conditions}, Transactions of the American Mathematical Society, 362(2)(2010), 1047-1096.

 \end{thebibliography}
\end{document}